\documentclass[11pt]{article}

\usepackage{fullpage} 

\usepackage[usenames,dvipsnames]{color} 
\usepackage{amsthm, amssymb, amsmath, amsfonts}
\usepackage{bm} 
\usepackage{bbm} 
\usepackage{tikz} 
	\usetikzlibrary{shapes,snakes} 
	\usetikzlibrary{arrows} 
		\newcommand{\midarrow}{\tikz \draw[-triangle 90] (0,0) -- +(.1,0);}
\usepackage{hyperref} 
\usepackage{stmaryrd} 

\theoremstyle{plain}
\newtheorem{theorem}{Theorem}[section]
\newtheorem{proposition}[theorem]{Proposition}
\newtheorem{corollary}[theorem]{Corollary}
\newtheorem{lemma}[theorem]{Lemma}
\theoremstyle{definition}
\newtheorem{definition}[theorem]{Definition}
\newtheorem{remark}[theorem]{Remark}
\newtheorem{example}[theorem]{Example}

\DeclareMathOperator{\il}{IL}
\DeclareMathOperator{\diag}{diag}
\DeclareMathOperator{\st}{SF}

\definecolor{oxfordblue}{rgb}{0.0, 0.13, 0.28}

\title{Inertias of Laplacian matrices of weighted signed graphs}
\author{K.~{Hassani Monfared}
\thanks{Department of Mathematics and Statistics, University of Victoria, Victoria, BC, Canada\ \ V8W 2Y2} \and 
		G.~MacGillivray$^{\ *\ \ddag}$
		 \and
		D.~D.~Olesky
		\thanks{Department of Computer Science, University of Victoria, Victoria, BC, Canada\ \ V8W 2Y2}
		 \and 
		P.~van den Driessche$^{\ *}$ \thanks{Research supported by NSERC Discovery Grant}}
\date{\today}		

\begin{document}
\maketitle

\begin{abstract}
	We study the sets of inertias achieved by Laplacian matrices of weighted signed graphs. First we characterize signed graphs with a unique Laplacian inertia. Then we show that there is a sufficiently small perturbation of the nonzero weights on the edges of any connected weighted signed graph so that all eigenvalues of its Laplacian matrix are simple. Next, we give upper bounds on the number of possible Laplacian inertias for signed graphs with a fixed flexibility $\tau$ (a combinatorial parameter of signed graphs), and show that these bounds are sharp for an infinite family of signed graphs. Finally, we provide upper bounds for the number of possible Laplacian inertias of signed graphs in terms of the number of vertices.
\end{abstract}



\bigskip\bigskip
\begin{center}
	\emph{\Large Dedicated to Charles R. Johnson}
\end{center}

\section{Introduction}
The signs of the real parts of the eigenvalues of a coefficient matrix in a system of linear ordinary differential equations determine the stability of the dynamical system that it is describing. 
In this paper, we focus on systems with a coefficient matrix that is a symmetric Laplacian matrix, 
for which all eigenvalues are real.
In order to understand the behaviour of such a dynamical system,  
the number of positive, negative, and zero eigenvalues of the coefficient matrix, 
together known as the inertia of the matrix, are studied. 
The sign pattern of the Laplacian matrix is described by a graph with positive and negative edges. 
The spectrum of the Laplacian of a signed graph reveals important information about the underlying dynamics, for example see \cite{ADKM08, CWLBQ17, PSM16}.

A \emph{signed graph} is a simple graph $G = (V(G), E(G))$ together with an assignment of the signs 
$+$ or $-$ to its edges.
An edge is \emph{positive} if it is assigned the sign $+$, and \emph{negative} if it is assigned the sign $-$.
We assume throughout that a signed graph $G$ has vertex set $V(G) = \{1, 2, \ldots, n\}$, and a total of $m$ 
edges of which $m_+$ are positive and $m_-$ are negative.
The number of components of $G$ is $c = c(G)$.
There has recently been much interest
in the Laplacian spectrum and adjacency spectrum of signed graphs; 
see, for example,  \cite{BB19, WWT18} and references therein.

We study the Laplacian inertia of a more general graph in which edges are signed and weighted, as in \cite{PSM16, BD14, PPP18, BDK15, CWLBJQ16}.
A \emph{weighted signed graph}, $\Gamma = (G, \gamma)$, 
is a signed graph $G$ with an assignment $\gamma$ of non-zero real numbers to the edges
of $G$ such that, for each edge $e$, 
the sign of the number $\gamma(e)$ is the same as the sign of $e$.
We call $\gamma$ a \emph{consistent weighting} of $G$.
The weight  assigned to the edge $e = ij$ is denoted by $\gamma_{ij}$  rather than $\gamma(e)$.

	\begin{definition}
		The \emph{weighted adjacency matrix} of the weighted signed graph $\Gamma = (G, \gamma)$ is the symmetric $n \times n$ matrix $A = \left[  a_{ij} \right]$ where
		\begin{equation*}
			a_{ij} = 
			\begin{cases} 
				\gamma_{ij} & \text{ if } \quad ij \in E(G)\\ 
				0 &  \text{ if } \quad ij \not\in E(G).
			\end{cases}
		\end{equation*}
The \emph{(weighted) degree} $d_i$ of vertex $i$ is the sum of the weights of the edges adjacent to it. 
		Let $D$ be the diagonal matrix $\diag(d_1, d_2, \ldots, d_n)$. 
		The \emph{weighted Laplacian} of $\Gamma$ is 
		\begin{equation*}
			L = L(\Gamma) = L(G,\gamma) = A - D.
		\end{equation*}
		
	\end{definition}
	
Let $\Gamma = (G, \gamma)$ be a weighted signed graph.
Then the matrix $L(\Gamma)$ is symmetric, has 0 row and column sums, 
has 0 as an eigenvalue with multiplicity at least one, and 
the all-ones-vector, $\mathbbm{1}$, is an eigenvector corresponding to the zero eigenvalue.
When $G$ has only positive edges, $L(\Gamma)$ is a symmetric negative semi-definite matrix where the multiplicity of the zero eigenvalue (nullity) is equal to the number of connected components of $G$.

\begin{definition}
Let $G$ be a signed graph.
The \emph{Laplacian inertia} of a weighted signed graph $\Gamma = (G, \gamma)$ on $n$ vertices is
$IL(\Gamma) = (n_+, n_-, n_0)$,
 where $n_+$ is the number of positive eigenvalues of $L(\Gamma)$,
$n_-$ is the number of negative eigenvalues of $L(\Gamma)$,
 and $n_0 = n - n_+ - n_- \geq1$
is the multiplicity of 0 as an eigenvalue of $L(\Gamma)$.
We say that $G$ \emph{achieves} the Laplacian inertia $(n_+, n_-, n_0)$ and that $\Gamma$ \emph{realizes} it.
\end{definition}

We are interested in the set of all possible Laplacian inertias achieved by a signed graph $G$,
that is, the set of all triples $(n_+, n_-, n_0)$ that can be realized by a weighted signed graph
$\Gamma = (G, \gamma)$, where $\gamma$ is a consistent weighting of $G$.
If $IL(\Gamma)$ is a fixed triple for all consistent weightings $\gamma$ of $G$, then
$G$ has \textit{unique Laplacian inertia}.

Most of the remaining definitions and notation of this section are similar to those of \cite{BD14},
but are stated for signed graphs that may be disconnected.
	
Let $G$ be a signed graph.
Let $E_+$ be the set of positive edges of $G$, and $E_-$ be the set of negative edges of $G$.
Note that $|E_+| = m_+$ and $|E_-| = m_-$.
We denote by $G_+$ the spanning subgraph of $G$ with edge set $E_+$, 
and by $G_-$  the spanning subgraph of $G$ with edge set $E_-$.
The number of components of $G_+$ and $G_-$ are denoted by $c_+ = c_+(G)$ and $c_- = c_-(G)$, 
respectively.

	Applying Theorem 2.10 of \cite{BD14} to each connected component of a weighted signed graph gives the following result.
	
	\begin{theorem}
	\label{cor_inertia_bounds_disconnected}
	Let $G$ be a signed graph.
	For any consistent weighting $\gamma$, the following inequalities hold for the 
	components of the Laplacian inertia $(n_+, n_-, n_0)$ of the weighted signed graph
	$\Gamma = (G, \gamma)$:
			\begin{align}
			\label{bound1} c_+ - c &\leq n_+ \leq n - c_-\\
			\label{bound2} c_- - c &\leq n_- \leq n - c_+\\
			\label{bound3} c 		&\leq n_0 \leq n + 2c - c_- - c_+.
		\end{align}
	\end{theorem}
	Note that the proof of the above theorem in \cite{BD14} uses a $1$-parameter family $\Gamma(t)$ introduced below, but it is clear that the results hold for any consistent weighting of the edges of $G$.
	
	\begin{definition}
		The \emph{flexibility} of a 
		 (weighted) signed graph on $n$ vertices is $\tau = n + c - c_+ - c_-$.
	\end{definition}

Let $\Gamma = (G, \gamma)$ be a weighted signed graph.
We denote by $\Gamma_+$ the weighted signed graph that consists of $G_+$ and the restriction of $\gamma$ to $E_+$, and 
by $\Gamma_-$ the weighted signed graph that consists of $G_-$ and the restriction of $\gamma$ to $E_-$.
For a real number $t$, define $\Gamma(t) = \Gamma_+ + t \Gamma_-$.
The eigenvalues $\lambda_1(t) = \lambda_2(t) = \cdots = \lambda_c(t) = 0 $, and $\lambda_{c+1}(t) \leq \lambda_{c+2}(t) \leq \cdots \leq \lambda_n(t)$ of $L(\Gamma(t))$ are  
differentiable functions of the parameter $t$.
Observe that the weighting of $\Gamma(t)$ is a consistent weighting of $G$ if only if $t > 0$.

Define  the \emph{crossing polynomial}
$$M(\Gamma(t)) = \dfrac{(-1)^{n-1}}{n} \prod\limits_{i=c+1}^{n} \lambda_i(t).$$
Up to a multiplicative factor, $M(\Gamma(t))$ is the coefficient of the 
 term with exponent $c$
in the characteristic polynomial of $L(\Gamma(t))$.
The following results follow from those of Bronski and Deville \cite{BD14} for connected (weighted) signed graphs, and describe properties of $M(\Gamma(t))$ and justify calling it the crossing polynomial.

For a weighted signed graph $\Gamma$, let
$SF_k(\Gamma)$ be the set of all maximal spanning forests (with $n-c$ edges) of $\Gamma$ having exactly $k$ edges in $\Gamma_-$, 
and $\pi(F)$ be the product of the edge weights of the weighted maximal spanning forest $F$ of $\Gamma$. 

\begin{theorem}[\cite{BD14}, Lemmas 2.14, 2.16, 2.18]
$ $\\
\begin{enumerate}
\item[(1)]  $M(\Gamma(t)) = \sum\limits_{k = c_+ - c}^{n - c_-} a_k (-t)^k$,
where $a_k = \sum\limits_{F \in \st_k(\Gamma)} |\pi(F)|.$
The first and last coefficients,  namely $a_{c_+-c}$ and $a_{n-c_-}$, are strictly positive.
\item[(2)]  The roots of the polynomial $M(\Gamma(t))$ are real and non-negative.
\item[(3)]  For $i \geq c+1$, the eigenvalues $\lambda_i(t)$ are non-decreasing functions of $t$ that
cross zero transversely, that is, if $\lambda_i(t) = 0$, then $\lambda^\prime_i(t) > 0$.
\end{enumerate}

\label{lem_two_fourteen}
\label{lem_two_sixteen}
\label{lem_two_eighteen}
\end{theorem}
			
As $t$ increases, some negative eigenvalues of $L(\Gamma(t))$ become zero 
where $M(\Gamma(t)) = 0$, and then become positive.
The (simple/multiple) zeros of $M(\Gamma(t))$ are called the (single/multiple) 
\emph{eigenvalue crossings} of $L(\Gamma(t))$. In particular, if a zero of $M(\Gamma(t))$ has multiplicity $k \geq 2$, then we call it a $k$-crossing.
Note that if either of $m_+$ or $m_-$ is zero, then $L(\Gamma(t))$ is either 
positive semidefinite or negative semidefinite, respectively, 
for all $t > 0$ and there are no crossings. 
In the following lemma we show that the flexibility, $\tau$, equals the number of eigenvalue crossings as $t$ increases.

	\begin{lemma}
	\label{lem_tau_is_crossing}
Let $\Gamma = (G, \gamma)$.
There are exactly $\tau = n + c - c_- - c_+$ positive values of $t$ such that 
$M(\Gamma(t)) = 0$, counting multiplicities. 
Thus, if all of these zeros are simple, then $\tau$ is the number of single crossings. 
Furthermore, $\tau = 0$ if and only if $G$ has the unique Laplacian inertia $(c_ + -c,\, c_ - - c,\, c)$.
	\end{lemma}
	\begin{proof}
	By Theorem \ref{lem_two_sixteen} (2) the zeros of $M(\Gamma(t))$ are all real and non-negative. 
	Also, since 
		\begin{equation*}
			M(\Gamma(t)) = \sum_{k = c_+ - c}^{n - c_-} a_k (-t)^k 
			= (-t)^{c_+ - c} \sum_{k=0}^{n + c - c_- - c_+} a_{k + c_+ - c} (-t)^k,
		\end{equation*}
by Theorem  \ref{lem_two_fourteen} (1) it has exactly $c_+ - c$ zeros equal to $0$. 
Thus, independent of the consistent weighting $\gamma$, the rest of the zeros are positive, and there are $\tau = n + c - c_- - c_+$ of them. Finally, there are no crossings if and only if $\tau = 0$ which, 
since equality then holds in the lower and upper bounds in Theorem \ref{cor_inertia_bounds_disconnected},
 is equivalent to $G$ having the unique Laplacian inertia $(c_ + -c,\, c_ - - c,\, c)$.
	\end{proof}

Let $r \in \mathbb{R}^+$, and let $\Gamma = (G,\gamma)$ be a weighted signed graph on $n$ vertices. 
Let $r\Gamma$ be obtained from $\Gamma$ by multiplying the weight of each edge by $r$. 
Also, let $r^- \Gamma$  be  obtained from $\Gamma$ by multiplying the weight of only the negative edges by $r$.
Note that 
$c_\pm(r^- \Gamma) = c_\pm(r\Gamma) = c_\pm(\Gamma) = c_\pm(G) = c_\pm$ 
because the number of components depends only on the signed graph $G$.
	
	\begin{proposition}
	$M(r\Gamma(t)) = r^{n-1} M(\Gamma(t))$ and $M(r^-\Gamma(t)) = M(\Gamma(rt))$.
	\end{proposition}
	\begin{proof}
	Let $M(\Gamma(t)) = \sum\limits_{k=c_+ - c}^{n - c_-} a_k (-t)^k$, and 
	$M(r\Gamma(t)) = \sum\limits_{k=c_+ - c}^{n - c_-} b_k (-t)^k$. 
	Then for each $k$,
		 \begin{equation*}
		 	b_k = \sum_{F \in \st_k(r\Gamma)} |\pi(F)| = \sum_{F \in \st_k(\Gamma)} r^{n-1} |\pi(F)| = r^{n-1} a_k.
		 \end{equation*}
		 Hence, $M(r\Gamma(t)) = r^{n-1} M(\Gamma(t))$.

	Let $M(r^-\Gamma(t)) = \sum\limits_{k=c_+ - c}^{n - c_-} u_k (-t)^k$.
	For each $k$
		 \begin{equation*}
		 	u_k = \sum_{F \in \st_k(r^- \Gamma)} |\pi(F)| = \sum_{F \in \st_k(\Gamma)} r^{k} |\pi(F)| = r^{k} a_k.
		 \end{equation*}
		 Hence, $M(r^-\Gamma(t)) = M(\Gamma(rt))$.

	\end{proof}
	
	\begin{example}
	\label{ex_kay_three_crossing}
		Let $\Gamma$ be the complete graph on $3$ vertices with one positive edge and two negative edges, consistently weighted $\pm 1$. Then $M(\Gamma(t)) = t(t - 2)$. 
		Also, for any $r > 0$, $M(r\Gamma(t)) = r^2 t(t - 2)$, 
		and  $M(r^-\Gamma(t)) = rt(rt - 2)$.
	\end{example}

\section{Signed Graphs with Unique Laplacian Inertias}
	
	In this section we characterize the signed graphs that have a unique Laplacian inertia.
 We first show that forests have a unique Laplacian inertia. Then we will argue that creating a cycle in the graph with all edges of the same sign will not increase the number of Laplacian inertias, but a cycle with edges of both signs will introduce new Laplacian inertias. Putting all this together, we prove that graphs with unique Laplacian inertias are precisely the ones that have no `mixed' blocks.

	\begin{theorem}
		\label{thm_forest_inertia}
		Let $\Gamma = (G, \gamma)$ be a weighted signed forest with $m_+$ positive edges and $m_-$ negative edges. Then 
		\begin{equation*}
			\il(\Gamma) = ( m_-, \, m_+, \, n - m_+ - m_-) = (c_+ - c, c_- - c, c).
		\end{equation*}
	\end{theorem}
	\begin{proof}
For a signed forest with $c$ components,
$n - c = m_+ + m_-$, and 	$c_\pm = m_\mp + c$.
Hence the flexibility is 
$$\tau = n + c - c_+ - c_-  = (c + m_+ + m_-) + c - (m_- + c) - (m_+ + c) = 0.$$
		Hence there are no eigenvalue crossings, and
the Laplacian inertia is as given in Lemma \ref{lem_tau_is_crossing}. 
	\end{proof}
	
	\begin{corollary}
	\label{lem_nullity_of_forests}
		Let $\Gamma = (G,\gamma)$ be a weighted signed forest with $c$ connected components. Then the nullity of the Laplacian matrix of $\Gamma$ is $c$.
	\end{corollary}

Let $G$ be a signed forest on $n$ vertices and with a total of $m$ edges and $c$ connected components. 
For any consistent weighting $\gamma$ of $G$, the Laplacian of the weighted signed graph $\Gamma = (G, \gamma)$
has nullity  $n_0 = c = n - m$. 
This is the same as the nullity of the (usual) Laplacian of any  graph with all positive edges and $c$ components.
Also, signed forests with the same number of vertices and the same number of positive and negative edges have the same Laplacian inertia.
This is due to the obvious fact that trees have unique spanning trees, which results in only one term in the crossing polynomial. 

Below, we first consider the crossing polynomial of a cycle, 
and then show that the presence of both positive and
negative edges in a cycle results in eigenvalue crossings and more than one Laplacian inertia.
	
	\begin{example}
	Consider the signed graph $G$ in Figure \ref{fig1} where the solid black lines are positive edges and the dashed red line is a negative edge, and $\Gamma = (G,\gamma)$ with weights $a,b,c > 0$ and $d < 0$. (Note that, in this example, $c$ is a weight.)
		
		\begin{figure}
		\begin{center}
			\begin{tikzpicture}
				\node[draw, circle, black] (1) at(0,0) {$1$};
				\node[draw, circle, black] (2) at(0,2) {$2$};
				\node[draw, circle, black] (3) at(2,2) {$3$};
				\node[draw, circle, black] (4) at(2,0) {$4$};
				\draw[-, ultra thick, black] (1) -- node [midway,left] {$a$} 
											  (2) -- node [midway,above] {$b$}
											  (3) -- node [midway,right] {$c$} (4);
				\draw[dashed, ultra thick, red] (1) -- node [midway,below] {$d$} (4);
			\end{tikzpicture}
		\end{center}
		\caption{A weighted 4-cycle with three positive edges and one negative edge. Solid black lines are positive edges and the dashed red line is the negative edge.}
		\label{fig1}
		\end{figure}
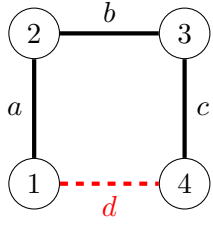
		
		By Theorem $\ref{lem_two_fourteen}$, the coefficients of $M(\Gamma(t))$ are $a_0 = abc$,
	$a_1 = -d (bc + ca + ab)$, $a_2 = a_3 = 0$.
		That is, 
		\begin{equation*}
			M(\Gamma(t)) = abc + d (bc + ca + ab) \; t,
		\end{equation*}
		with its only positive zero being
		\begin{equation*}
			t = \frac{abc}{-d (bc + ca + ab)} > 0.
		\end{equation*}
		The signed graph $G$ achieves Laplacian inertias $(0,3,1)$, $(0,2,2)$, and $(1,2,1)$ for any fixed positive values $a,b,c$, and varying negative values of $d$.
		\end{example}
		
A similar argument gives the following result.			
		
		\begin{theorem}
		Any signed cycle $C$ on $n$ vertices with $m_+ > 0$ positive edges and $m_- > 0$ negative edges has $c_+ = m_-$, $c_- = m_+$, $n = m_+ + m_-$, and $\tau = n + 1 - c_+ - c_- = 1$. 
		The polynomial $M(\Gamma(t))$ with $\Gamma = (C,\gamma)$ has only two terms with one positive zero, giving one eigenvalue crossing. That is,
		\begin{equation*}
			M(\Gamma(t)) = (-t)^{m_- - 1} (a_{m_- - 1}  - a_{m_-} t),
		\end{equation*}
		with a unique positive zero $t = \frac{a_{m_- - 1}}{a_{m_-}}$.
		\end{theorem}	

	\begin{theorem}
		The flexibility $\tau$ of a graph $\Gamma$ is zero if and only if 
		there is no cycle with both a positive edge and a negative edge.
		\label{thm_mixed_cycles}
	\end{theorem}
	\begin{proof}
	There is no cycle with both a positive edge and a negative edge if and only if
 every spanning tree has the same number of negative edges. 
 
	($\Leftarrow$) 
If every spanning tree of $\Gamma$ has exactly $k$ negative edges, then $a_{j} \neq 0$ if and only if $j = k$. Then $M(\Gamma(t)) = a_k (-t)^k$, i.e. it has only one term, giving $\tau = 0$. 
	
	($\Rightarrow$) Consider the contrapositive.
Assume there is a cycle $C$ with both a positive edge and a negative edge. 
Then there are spanning trees with a different number of negative edges.
Therefore $M(\Gamma(t))$ with $\Gamma = (C,\gamma)$ has more than one term, 
so by Lemma \ref{lem_tau_is_crossing},
it has a positive zero. Thus $\tau > 0$.
	\end{proof}
	
A \emph{block} of a graph $G$ is a maximal connected subgraph with no cut vertex.  
By Menger's Theorem (see \cite{Bondy}, Corollary 3.3.2)  two edges of a graph belong to a common cycle if and only if
they are in the same block.
Thus, a signed graph has a cycle with both a positive edge and a negative edge 
if and only if it has a block with both a positive edge and a negative edge.	
Theorem \ref{thm_mixed_cycles} can therefore be restated as follows.

	\begin{corollary}
	\label{cor_unique_iniertia_no_mixed_block}
		A signed graph has a unique Laplacian inertia if and only if no  block has 
		both a positive edge and a negative edge.
	\end{corollary}
	
	Since there is an efficient algorithm to find the blocks of a graph \cite{T72}, the
previous corollary gives an algorithm to test if a signed graph has a unique Laplacian inertia:
first find the blocks and then scan each one to test if all its edges have the same sign.
		
	\begin{theorem}
		Let $\Gamma = (G, \gamma)$ be a weighted signed graph with $c$ connected components and no block having both positive and negative edges. Assume $G$ has $k$ blocks of sizes $n^+_1, n^+_2,\ldots, n^+_k$ with only positive edges and $\ell$ blocks of sizes $n^-_1, n^-_2,\ldots, n^-_\ell$ with only negative edges. Then $G$ has the unique Laplacian inertia $(n_+, n_-, n_0)$ if and only if 
		\begin{align*}
			n_+ = \sum_{i=1}^\ell (n^-_i - 1) - c, \,
			n_- = \sum_{i=1}^k (n^+_i - 1) - c, \,
			\text{and } n_0 = c.
		\end{align*}
	\end{theorem}
	\begin{proof}
		By Corollary \ref{cor_unique_iniertia_no_mixed_block}, $G$ has unique Laplacian inertia. Then by Lemma \ref{lem_tau_is_crossing}, $\tau = 0$ and the unique Laplacian inertia is $(c_+ - c, c_- - c, c)$.		Note that for $G$ with blocks as specified, $c_+ = \sum_{i=1}^\ell (n^-_i - 1) \geq c$ and $c_- = \sum_{i=1}^k (n^+_i - 1) \geq c$. 
	\end{proof}
	
	\begin{theorem}
		There is a signed graph $G$ with unique Laplacian inertia $(n_+,n_-,n_0)$ if and only if $n_\pm \geq 0$ and $n_0 \geq 1$. Furthermore, for any such graph $n_+ = c_+ - c$, $n_- = c_- - c$, and $n_0 = c$. 
	\end{theorem}
	\begin{proof}
		If $G$ has unique Laplacian inertia, then by Lemma \ref{lem_tau_is_crossing}, the unique Laplacian inertia is $(c_+ - c, c_- - c, c)$, which implies $n_\pm \geq 0$ and $n_0 \geq 1$. Conversely, given $(n_+, n_-, n_0)$, first realize $(n_+, n_-, 1)$ for a forest using Theorem $\ref{thm_forest_inertia}$, and then add $c-1$ isolated vertices to give $G$.
	\end{proof}

\section{Upper bounds for the number of Laplacian Inertias}

	It is easy to see that the set of (real symmetric) matrices with all eigenvalues simple is dense in the set of (real symmetric) matrices. In other words, there is a sufficiently small (symmetric) perturbation of any (real symmetric) matrix so that the perturbed matrix has all eigenvalues simple, that is, multiple eigenvalues are rare. In this section we first show that in the same sense multiple eigenvalues are rare for Laplacian matrices of weighted signed graphs. That is, there is a sufficiently small perturbation of the weights of a weighted signed graph resulting in only simple eigenvalues for its Laplacian. Then we give upper bounds on the number of Laplacian inertias for signed graphs with a fixed flexibility, and separately for signed graphs on a fixed number of vertices. Throughout, we discuss if these upper bounds are sharp.

\subsection{Perturbing a Laplacian matrix to obtain simple eigenvalues}
		Poignard et al.~in \cite{PPP18} point out that the simplicity of eigenvalues of the Laplacian matrices of weighted graphs plays an important role in guaranteeing ``good properties of the underlying dynamics such as exponentially and uniformly fast convergence towards synchronization in diffusively coupled networks and convergence to the stationary measure in random walks.'' In that paper they show that there is a sufficiently small perturbation of the nonzero weights of the edges of any connected weighted graph with all edges positive so that all eigenvalues of its Laplacian are simple. Here we show that their proof can be modified to prove that the same conclusion holds true for the Laplacian matrices of connected weighted \emph{signed} graphs. This modification is due to the signs of the edges, which imply that the number of connected components in a signed graph is no longer necessarily equal to the nullity of the Laplacian matrix of the graph. However, in the proof we only need to consider the nullity of the Laplacian of spanning signed forests of the graph, and Corollary \ref{lem_nullity_of_forests} guarantees that the nullity of the Laplacian of a weighted signed forest is still equal to its number of connected components.
		
		We refrain from providing a complete proof here, as it would mostly be a repetition of the proof of Theorem $3.1$ in \cite{PPP18}, and we encourage the interested reader to refer to the original paper. Instead we give an outline of their proof and what needs to be changed so that it works for connected weighted signed graphs. In what follows, for any consistent weighting $\gamma$ of a signed graph $G$ let $||\gamma|| = ||L(G,\gamma)||_2$, though any norm can be used.
		
		\begin{theorem}
			For any connected weighted signed graph $\Gamma = (G,\gamma)$ and a real number $\varepsilon > 0$ there is a consistent weighting $\gamma'$  of the edges of $G$ such that $||\gamma - \gamma'|| < \varepsilon$ and $L(G,\gamma')$ has only simple eigenvalues. In other words, the set of consistent weightings of a connected signed graph $G$ for which the Laplacian matrices have all of their eigenvalues simple is dense in the set of all consistent weightings of $G$.
		\end{theorem}
		\begin{proof}[Sketch of the proof]
			The proof has two main steps. First, given $\Gamma = (G, \gamma)$ obtain a consistent weighting $\gamma_1$ of $G$ with $||\gamma_1|| < \varepsilon$ so that $L(G,\gamma_1)$ has all eigenvalues simple. Second, show that for $\gamma' = \gamma + \frac{b}{a}\gamma_1$ (for some $a,b > 0$ and $\frac{b}{a} < 1$) the eigenvalues of $L(G,\gamma')$ are also simple. Note that since $||\gamma_1|| < \varepsilon$, then $||\gamma - \gamma'|| < \varepsilon$ and $\gamma'$ is also a consistent weighting of $G$.
			
			\begin{itemize}
				\item[Step 1.] Assume that $G$ has $n$ vertices and a longest path $P$ of $G$ is of length $k$. In order to obtain the consistent weighting $\gamma_1$, as described above, start with the path $P$ in $G$ and choose the weights of edges of $P$ sufficiently small and consistent with the signs in $G$. First, note that without loss of generality, the Laplacian matrix of $P$ is irreducible, symmetric, and tridiagonal. Hence all of its eigenvalues are simple. Now add the remaining $n-k$ vertices of $G$ to $P$ as isolated vertices. By Corollary \ref{lem_nullity_of_forests} the Laplacian matrix now has nullity $n-k+1$.
				
				Repetitively connect one of the remaining isolated vertices to the tree $P$ as a pendent vertex by choosing its weight sufficiently small and consistent with the signs in $G$ until a signed spanning tree $T$ of $G$ of arbitrarily small norm is obtained. In each step, adding an edge as described reduces the number of connected components by 1, and by Corollary \ref{lem_nullity_of_forests} the nullity decreases by $1$. It can be shown that the weight of the edge being added can be chosen sufficiently small so that all the nonzero eigenvalues of the corresponding Laplacian matrix remain simple. That is, at the end of this process, $T$ is a signed spanning tree of $G$ with a consistent weighting of arbitrarily small norm whose Laplacian has all eigenvalues simple. 
								
				Finally, all the edges of $G$ that are not yet in $T$ can be added to $T$ with sufficiently small weights to create a consistent weighting $\gamma_1$ of $G$ so that $L(G, \gamma_1)$ has all eigenvalues simple and $||\gamma_1|| < \varepsilon$. Note that the only difference in the sketch of the proof provided here from the proof of \cite{PPP18} is that in their case, the graph is not signed and all the eigenvalues are nonnegative, but here some of the eigenvalues might be negative.
				
				\item[Step 2.] The weighting $\gamma_1$ above is a consistent weighting of $G$ with norm less than $\varepsilon$ and $L_1 = L(G,\gamma_1)$ has all eigenvalues simple. It is clear that for positive 
				$\tilde{a}$ and sufficiently small positive $\tilde{b}$ a weighting $\tilde{\gamma} = \tilde{a} \gamma +  \tilde{b} \gamma_1$ 				
is also a consistent weighting of $G$. 
As in the proof of Theorem $3.1$ in \cite{PPP18}, it follows that the Laplacian matrix $\tilde{L} = L(G,\tilde{\gamma})$ has all eigenvalues simple. This can be shown by considering the discriminant function 
				\begin{center}
				\begin{tabular}{ccrcl}
					${\rm D}_{L,L_1}(\tilde{a},\tilde{b})$ &: 	&$\mathbb{R}^2$ 	&$\to$ 		&$\mathbb{R}$, \\
											  &		&$(\tilde{a},\tilde{b})$ 		&$\mapsto$	& $c(L,L_1) \displaystyle\prod_{i<j} (\lambda_i - \lambda_j) ^2$,
				\end{tabular}
				\end{center}
				where $c(L,L_1)$ is a constant and the $\lambda_i$ are the eigenvalues of the Laplacian matrix $\tilde{L}$. Note that ${\rm D}_{L,L_1}$ is a polynomial in the entries of $\tilde{L}$ and it is zero if and only if $\tilde{L}$ admits at least one multiple eigenvalue. Moreover, it is not identically zero over $\mathbb{R}^2$ since, for example, $(\tilde{a},\tilde{b}) = (0,1)$ gives $\tilde{\gamma} = \gamma_1$ and $L_1$ has all eigenvalues simple. Hence, 
there are  $a,b > 0$ such that $\displaystyle\frac{b}{a} < 1
$, and $L + \displaystyle\frac{b}{a} L_1$ has all eigenvalues simple. Choose $\gamma' =  \gamma + \displaystyle\frac{b}{a} \gamma_1$. Then $L' = L(G, \gamma')$ has all eigenvalues simple  and $||L' - L||_2 < \varepsilon$.
			\end{itemize}
		\end{proof}

	\subsection{Graphs with a fixed flexibility}
		So far we have shown that there is a sufficiently small perturbation of a given weighting of a connected weighted signed graph so that all eigenvalues of its Laplacian matrix are simple. This means that not only for graphs with unique Laplacian inertia is the nullity of the Laplacian equal to the number of connected components, but also, any weighting of a graph with $c$ connected components with the nullity of the Laplacian more than $c$ is a small perturbation away from a weighting that yields $n_0 = c$. In other words, the lower bound given by \eqref{bound3} is always realizable. As shown in \cite{BD14}, this lower bound is realized with the lower (resp. upper) and upper (resp. lower) bounds of \eqref{bound1} and \eqref{bound2}. 
		What we do not know yet is what other Laplacian inertias given by the bounds of \eqref{bound1}--\eqref{bound3} are realizable. Note that the difference between the upper and lower bounds in all three inequalities is equal to the flexibility $\tau$ of the graph. Here we introduce a lattice to visualize these bounds and to easily count how many different Laplacian inertias satisfy them, namely ${\tau + 2 \choose 2}$. We then show that for any given $\tau$ there is a connected signed graph that achieves all ${\tau + 2 \choose 2}$ possible Laplacian inertias.

		\begin{lemma}
			A connected signed graph with flexibility $\tau$ has at most ${\tau + 2 \choose 2}$ Laplacian inertias.
		\end{lemma}
		\begin{proof}
			Consider the plane with horizontal axis $n_+$ and vertical axis $n_-$. For a fixed $n$, each Laplacian inertia $(n_+, n_-, n_0)$ is represented by a point of a lattice in this plane, since $n_0 = n - n_+ - n_-$. Now for a given signed graph with flexibility $\tau$, assume all single and multiple crossings are realizable. One of the extreme Laplacian inertias given by Theorem \ref{cor_inertia_bounds_disconnected}, for example with $t > 0$ sufficiently small, has $n_- = n - c_+$ and $n_+ = c_+ - 1$. Each single crossing can be represented by going one step down and then one step right starting from this pair $(n_+, n_-)$. A $2$-crossing can be represented by going two steps down and then two steps right, and so on. Finally, a $\tau$-crossing can be represented by going $\tau$ steps down and then $\tau$ steps right. This creates a lattice in the $(n_+, n_-)$-plane with one point in the first row, two points in the second row, and so on. Finally, the lattice has $\tau+1$ points in the last row. Recall each point on the lattice represents a possible Laplacian inertia for the signed graph. Thus the number of Laplacian inertias for a given signed graph with flexibility $\tau$ is no more than the total number of the points on the lattice which is $1 + 2 + \cdots + (\tau+1) = {\tau+2 \choose 2}$. 
		\end{proof}
		
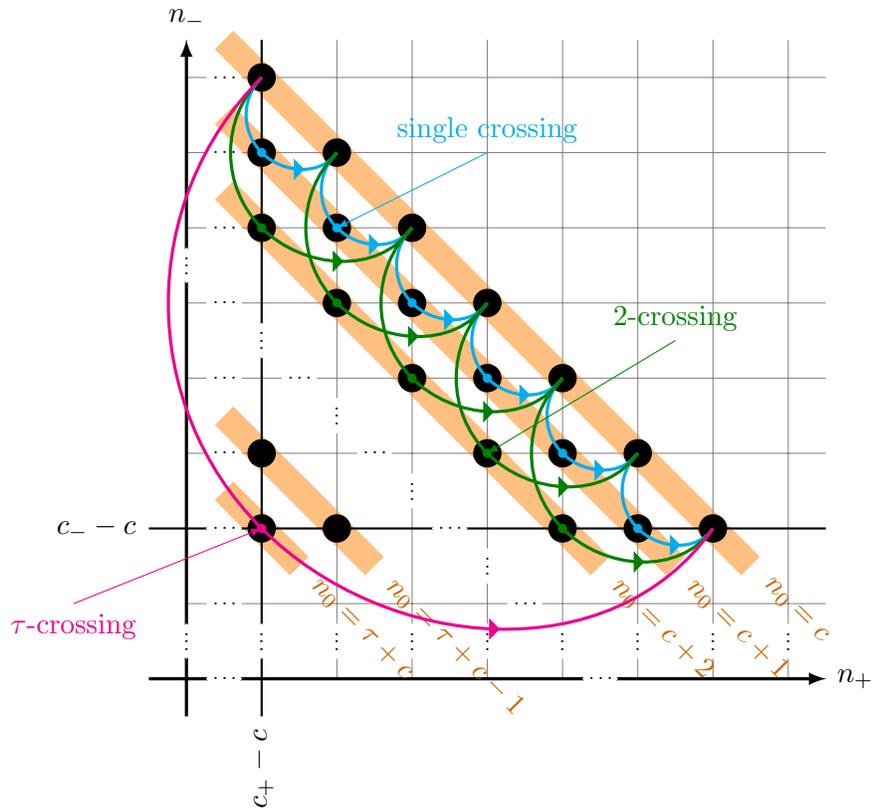
\begin{figure}			
			\begin{center}
				\begin{tikzpicture}[scale=1,
									nullity_text/.style={right, 
														 rotate=-45, 
														 scale=1,
														 opacity=1, 
														 orange!80!black}, 
									nullity_line/.style={line width=3.5mm, 
														 orange, 
														 opacity=0.5},
									single_crossing_node/.style={draw, 
																  circle, 
																  fill,
																  scale = .3,
																  cyan},
									double_crossing_node/.style={draw, 
																  circle, 
																  fill,
																  scale = .3,
																  green!50!black}]
						\draw[help lines,step=1] (0,0) grid (8.5,8.5);					
						\draw[-latex, very thick] (-.5,0) -- (8.5,0) node[right] {$n_+$};
						\draw[-latex, very thick] (0,-.5) -- (0,8.5) node[above] {$n_-$};					
					
						\draw [nullity_line] (0.5,8.5) -- (7.5,1.5) node [nullity_text] {$n_0 = c$};;		
						\draw [nullity_line] (0.5,7.5) -- (6.5,1.5) node [nullity_text] {$n_0 = c+1$};;	
						\draw [nullity_line] (0.5,6.5) -- (5.5,1.5) node [nullity_text] {$n_0 = c+2$};;	
						\draw [nullity_line] (0.5,3.5) -- (2.5,1.5) node [nullity_text] {$n_0 = \tau+c-1$};;	
						\draw [nullity_line] (0.5,2.5) -- (1.5,1.5) node [nullity_text] {$n_0 = \tau+c$};;		
						\draw[-, thick] (-.5,2) -- (8.5,2);
						\node[] () at (-1.2,2) {$c_- - c$};
						
						\draw[-, thick] (1,-.5) -- (1,8.5);
						\node[rotate=90, anchor=center] () at (1,-1.2) {$c_+ - c$};
					\foreach \y in {0,1,...,8}
					{
						\node[draw, diamond,fill,white,scale=.9] () at (.5,\y) {};
						\node[scale = .7] () at (.53,\y) {$\cdots$};
						
						\node[draw, diamond,fill,white,scale=.9] () at (\y,.5) {};
						\node[scale = .7] () at (\y,.57) {$\vdots$};
					}
						\pgfmathsetmacro{\h}{4.5}
						\pgfmathsetmacro{\z}{9-\h}
						\foreach \y in {-1,0,1,...,\h}
						{
							\node[draw, diamond,fill,white,scale=.9] () at (1+\y,\h-\y) {};
							\node[scale = .7] () at (1+\y,.07+\h-\y) {$\vdots$};
						}
						\foreach \y in {-.5,.5,...,\h}
						{
							\node[draw, diamond,fill,white,scale=.9] () at (1+\y,\h-\y) {};
							\node[scale = .7] () at (1.03+\y,\h-\y) {$\cdots$};
						}
						\foreach \y in {0,1,...,8}
						{
							\node[draw, diamond,fill,white,scale=.9] () at (.5,\y) {};
							\node[scale = .7] () at (.53,\y) {$\cdots$};
						}
						\foreach \h in {1,2,5  ,6,7}
						{
							\pgfmathsetmacro{\z}{9-\h}
							\pgfmathsetmacro{\g}{\h-1}
							\foreach \y in {0,...,\g}
							{
								\node[draw, circle, fill] (start) at (1+\y,1+\h-\y) {};
							}
						}
					\foreach \x in {7,6,...,2}
					{
						\node[single_crossing_node] () at (\x-1,9-\x) {};
						\draw[-, very thick, cyan] (\x,9-\x) 
													to[out=235,in = -45] 
													node {\midarrow} 
													(\x-1,9-\x) 
													to[out=135,in = -135] 
													(\x-1,10-\x);
					}
					\draw[latex-, cyan] (2,6) -- (4,7) node[above] {single crossing};
					
					\foreach \x in {7,6,...,3}
					{
						\node[double_crossing_node] () at (\x-2,9-\x) {};
						\draw[-, very thick, green!50!black] (\x,9-\x) 
												 to[out=235,in = -45] 
												 node {\midarrow} 
												 (\x-2,9-\x) 
												 to[out=135,in = -135] 
												 (\x-2,11-\x);
					}
					\draw[latex-, green!50!black] (4,3) -- (6.5,4.5) node[above] {$2$-crossing};
					
					\node[draw, 
						  circle, 
						  fill,
						  scale = .3,
						  magenta] () at (1,2) {};
					\draw[-, very thick, magenta] (7,2) 
												 to[out=235,in = -45] 
												 node {\midarrow} 
												 (1,2) 
												 to[out=135,in = -135] 
												 (1,8);
					\draw[latex-, magenta] (1,2) -- (-1.5,1) node[below] {$\tau$-crossing};
				\end{tikzpicture}
			\end{center}
		\caption{The lattice for a connected signed graph $G$ on $n$ vertices with flexibility $\tau$, $c_+$ connected components of $G_+$ and $c_-$ connected components of $G_-$.}		
		\label{fig2}
		\end{figure}

		The lattice for a connected signed graph $G$ on $n$ vertices with flexibility $\tau$, $c_+$ connected components of $G_+$ and $c_-$ connected components of $G_-$ is shown in Figure \ref{fig2}.
	
		\begin{definition}
			For two weighted signed graphs $\Gamma_1 = (G_1, \gamma_1)$ and $\Gamma_2 = (G_2, \gamma_2)$,  define $\Gamma_1 \cdot \Gamma_2$ to be any weighted signed graph obtained by identifying one vertex of $G_1$ with one vertex of $G_2$ and keeping the weights of the corresponding edges.
		\end{definition}
	
		\begin{theorem}\label{thm_crossing_dot}
			For two weighted signed graphs $\Gamma_1 = (G_1, \gamma_1)$ and $\Gamma_2 = (G_2, \gamma_2)$,
			\begin{equation*}
				M(\Gamma_1 \cdot \Gamma_2) = M(\Gamma_1) M(\Gamma_2).
			\end{equation*}
		\end{theorem}
		\begin{proof}
			Let
			\begin{align*}
				M(\Gamma_1 \cdot \Gamma_2) = \sum a_k (-t)^k, \quad
				M(\Gamma_1) = \sum b_k (-t)^k, \text{ and } 
				M(\Gamma_2) = \sum c_k (-t)^k.
			\end{align*}
			We want to show for all $k$ that
			\begin{equation*}
				a_k  = \sum_{r=0}^{k} b_r c_{k-r}.
			\end{equation*}
			For any fixed $k$ assume that $F \in \st_k(\Gamma_1 \cdot \Gamma_2)$, and let $F_1$ and $F_2$ be restrictions of $F$ to  the vertices of $\Gamma_1$ and $\Gamma_2$, respectively. Then $F_1 \in \st_r(\Gamma_1)$ and $F_2 \in \st_{k-r}(\Gamma_2)$ for some $0 \leq r \leq k$. Hence,
			\begin{equation*}
				\pi(F) = \pi(F_1) \pi(F_2),
			\end{equation*}
			and thus 
			\begin{align*}
				a_k &= \sum_{F \in \st_k(\Gamma)} | \pi(F) |\\
					&= \sum_{r=0}^{k} \sum_{\substack{F_1 \in \st_r(\Gamma_1) \\ F_2 \in \st_{k-r}(\Gamma_2)}} | \pi(F_1) \pi(F_2) | \\
					&= \sum_{r=0}^{k} \left( \sum_{F_1 \in \st_r(\Gamma_1)} |\pi(F_1)| \right) \left( \sum_{F_2 \in \st_{k-r}(\Gamma_2)} | \pi(F_2) | \right) \\
					&= \sum_{r=0}^{k} b_r c_{k-r}.
			\end{align*}
		\end{proof}

		\begin{example}
		\label{ex_k3_dots}
			Let $G$ be a $K_3$ with two negative edges and one positive edge as shown in Figure \ref{fig3}(A), and $\Gamma = (G,\gamma)$ have consistent weights $\pm 1$. The flexibility of $G \cdot G$ is $\tau = 2$ with $0 \leq n_- \leq 2$, $2 \leq n_+ \leq 4$, and the number of Laplacian inertias for $G$ is no more than ${\tau+2 \choose 2} = 6$, of which $5$ of them with $n_0 \leq 2$ are guaranteed by single crossings, and the one with $n_0 = 3$ can be obtained by weights $\gamma$ of $\Gamma$ as follows. By Example $\ref{ex_kay_three_crossing}$
			\begin{align*}
				M(\Gamma(t)) &= t(t-2),\\
				M(2^-\Gamma(t)) &= 2t(2t - 2), \text{ and } \\
				M(0.5^-\Gamma(t)) &= 0.5 t(0.5 t - 2).
			\end{align*}
			Thus, by Theorem $\ref{thm_crossing_dot}$
			\begin{equation*}
				M((2^-\Gamma \cdot \Gamma) (t)) = M(2^-\Gamma(t)) M(\Gamma(t)) = 2t^2(2t - 2)(t-2).
			\end{equation*}
			As $t$ increases from a small positive number to infinity, the weighted signed graph $(2^-\Gamma \cdot \Gamma) (t)$ realizes five different Laplacian inertias with $n_0 = 1, 2$. Also
			\begin{equation*}
				M((\Gamma \cdot \Gamma) (t)) = M(\Gamma(t)) M(\Gamma(t)) = t^2(t - 2)^2.
			\end{equation*}
			Thus $\Gamma \cdot \Gamma$ has a $2$-crossing at $t = 2$ as shown in Figure \ref{fig3}(B).
			\begin{figure}
			\begin{center}
				\begin{tikzpicture}
				\node[] () at (-.5,3) {(A)};
				\node[] () at (5.5,3) {(B)};
					\node[draw, circle, black] (1) at(0,0) {$1$};
					\node[draw, circle, black] (2) at(0,2) {$2$};
					\node[draw, circle, black] (3) at(2,0) {$3$};
					\node[draw, circle, black] (4) at(2,2) {$4$};
					\node[draw, circle, black] (5) at(4,0) {$5$};
					\draw[-, ultra thick, black] (1) -- (2);
					\draw[-, ultra thick, black] (3) -- (4);
					\draw[dashed, ultra thick, red] (1) -- (3) -- (2);
					\draw[dashed, ultra thick, red] (3) -- (5) -- (4);
					
					\node[scale=.3] () at (9	,1) {\includegraphics[]{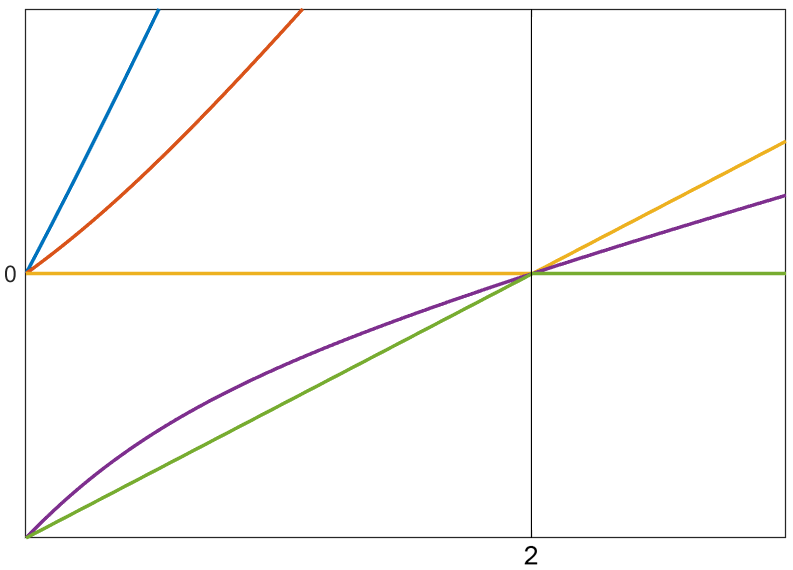}};
				\end{tikzpicture}
			\end{center}
			\caption{(A) Signed graph $G \cdot G$ where $G = K_3$ with one positive edge and two negative edges. (B) A $2$-crossing of a consistent weighting of $G \cdot G$ with weights $\pm 1$.}
			\label{fig3}
			\end{figure}
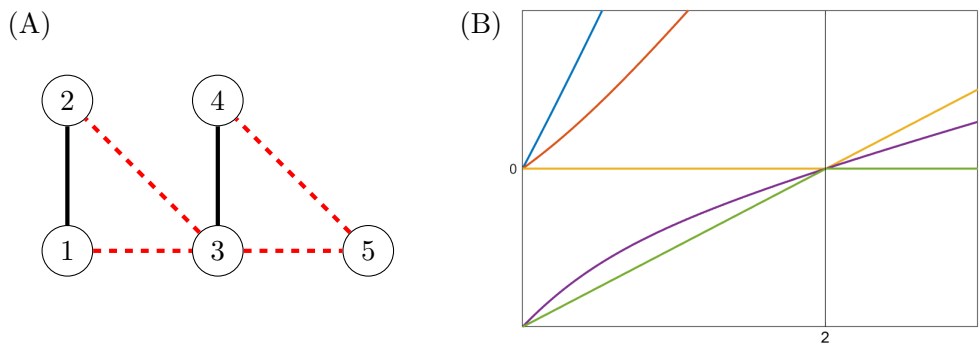
			The $6$ Laplacian inertias achieved by $G \cdot G$ are $(2,2,1)$, $(3,1,1)$, $(4,0,1)$, 
							    $(2,1,2)$, $(3,0,2)$, and
							    $(2,0,3)$.
			
			Similarly, the flexibility of $G \cdot G \cdot G$ (shown in Figure \ref{fig4}(A)) is $\tau = 3$ with $0 \leq n_- \leq 3$, $3 \leq n_+ \leq 6$,  and the number of Laplacian inertias for $G \cdot G \cdot G$ is at most ${\tau+2 \choose 2} = 10$, of which $7$ of them with $n_0 \leq 2$ are guaranteed by single crossings, the one with $n_0 = 4$ can be obtained by using the edge weights of $\Gamma$ since
			\begin{equation*}
				M((\Gamma \cdot \Gamma \cdot \Gamma) (t)) = t^3(t - 2)^3.
			\end{equation*}
			Hence, $\Gamma \cdot \Gamma \cdot \Gamma$ has a $3$-crossing at $t = 2$ as shown in Figure \ref{fig4}(B).
			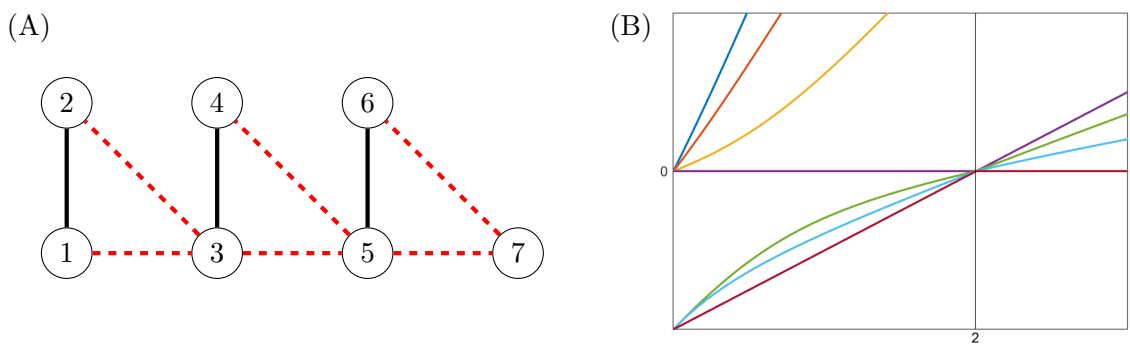
\begin{figure}
			\begin{center}
				\begin{tikzpicture}
				\node[] () at (-.5,3) {(A)};
				\node[] () at (7.5,3) {(B)};
					\node[draw, circle, black] (1) at(0,0) {$1$};
					\node[draw, circle, black] (2) at(0,2) {$2$};
					\node[draw, circle, black] (3) at(2,0) {$3$};
					\node[draw, circle, black] (4) at(2,2) {$4$};
					\node[draw, circle, black] (5) at(4,0) {$5$};
					\node[draw, circle, black] (6) at(4,2) {$6$};
					\node[draw, circle, black] (7) at(6,0) {$7$};
					\draw[-, ultra thick, black] (1) -- (2);
					\draw[-, ultra thick, black] (3) -- (4);
					\draw[-, ultra thick, black] (5) -- (6);
					\draw[dashed, ultra thick, red] (1) -- (3) -- (2);
					\draw[dashed, ultra thick, red] (3) -- (5) -- (4);
					\draw[dashed, ultra thick, red] (5) -- (7) -- (6);
					
					\node[scale=.3] () at (11,1) {\includegraphics[]{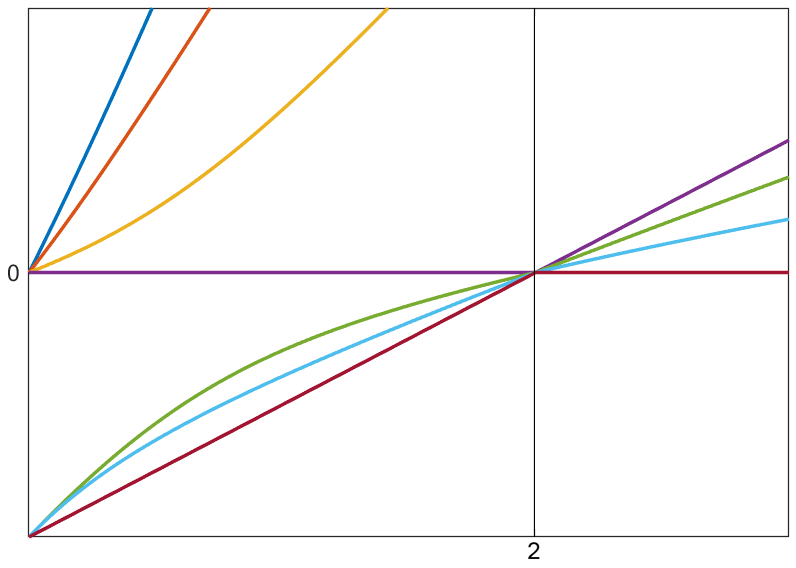}};
				\end{tikzpicture}
			\end{center}
			\caption{(A) Signed graph $G \cdot G \cdot G$ where $G = K_3$ with one positive edge and two negative edges. (B) A $3$-crossing of a consistent weighting of $G \cdot G \cdot G$ with weights $\pm 1$.}
			\label{fig4}
			\end{figure}
			Finally, the remaining two Laplacian inertias with $n_0 = 3$ can be obtained by
			\begin{equation*}
				M((2^-\Gamma \cdot \Gamma \cdot \Gamma) (t)) = 4t^3(t - 1)(t - 2)^2,
			\end{equation*}
			and 
			\begin{equation*}
				M((0.5^-\Gamma \cdot \Gamma \cdot \Gamma) (t)) = 0.25 t^3(t - 4)(t - 2)^2.
			\end{equation*}
			Thus, $2^-\Gamma \cdot \Gamma \cdot \Gamma$ and $0.5^-\Gamma \cdot \Gamma \cdot \Gamma $ have a single crossing (at $t = 1$ and $t = 4$, respectively), and both have a $2$-crossing at $t = 2$ as shown in Figure \ref{fig5}.
			\begin{figure}
			\begin{center}
				\includegraphics[width = .40\linewidth]{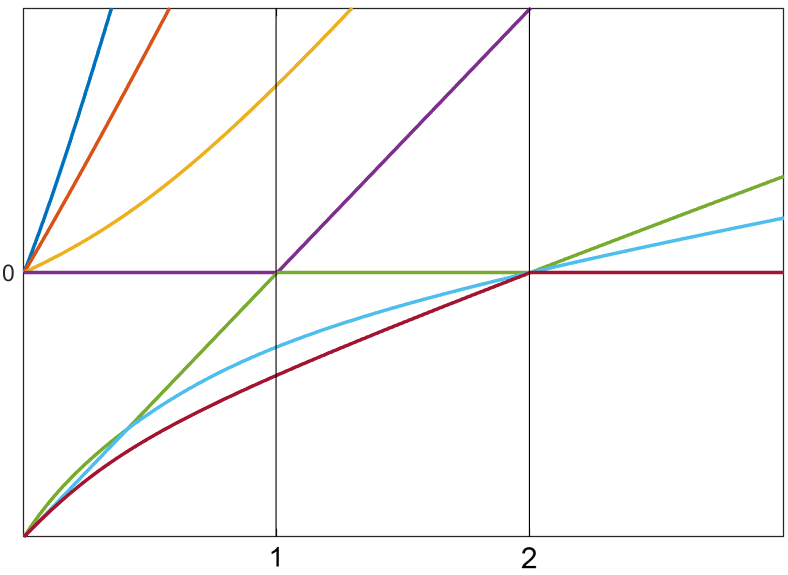} \hspace{1cm}
				\includegraphics[width = .40\linewidth]{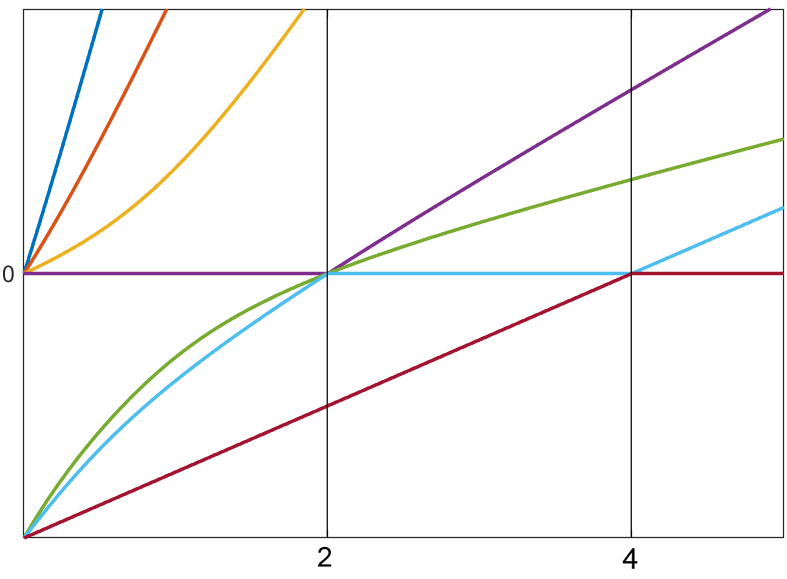}
			\end{center}
			\caption{Two different $2$-crossings of two particular consistent weightings of $G \cdot G \cdot G$ with $G$ as in Figure \ref{fig4}.}
			\label{fig5}
			\end{figure}
			The $10$ Laplacian inertias achieved by $G \cdot G \cdot G$ are  $(3,3,1)$, $(4,2,1)$, $(5,1,1)$, $(6,0,1)$, $(3,2,2)$, $(4,1,2)$, $(5,0,2)$, $(3,1,3)$, $(4,0,3)$, and $(3,0,4)$.
		\end{example}
		
		\begin{remark}
		\label{cor_mink_sum_inertias}
			For any two weighted signed graphs $\Gamma_1 = (G_1, \gamma_1)$ and $\Gamma_2 = (G_2, \gamma_2)$, the set of all possible pairs $(n_+, n_-)$ in the Laplacian inertias of $\Gamma_1 \cdot \Gamma_2$ is a subset of the Minkowski sum (the sumset) of the sets of possible pairs $(n_+, n_-)$ in the Laplacian inertias of $\Gamma_1$ and $\Gamma_2$. 
		\end{remark}
	
		\begin{theorem}
		\label{thm_tau_bound}
		For each $\tau \geq 0$ there is a signed graph with flexibility $\tau$ that achieves ${\tau+2 \choose 2}$ Laplacian inertias.
		\end{theorem}
		\begin{proof} 
			Let $G$ be a $K_3$ with two negative edges and one positive edge. Let $\Gamma = (G,\gamma)$ with all consistent weights $\pm 1$. Then the signed graph $G \cdot G \cdot \cdots \cdot G$, where $G$ is repeated $k$ times, has $n = 2k+1$ vertices with $c_- = 1$ and $c_+ = k+1$. Hence $\tau = k$. The pairs $(n_+, n_-)$ of the set of Laplacian inertias achieved by this signed graph are given in the lattice in Figure \ref{fig6}.	
			
			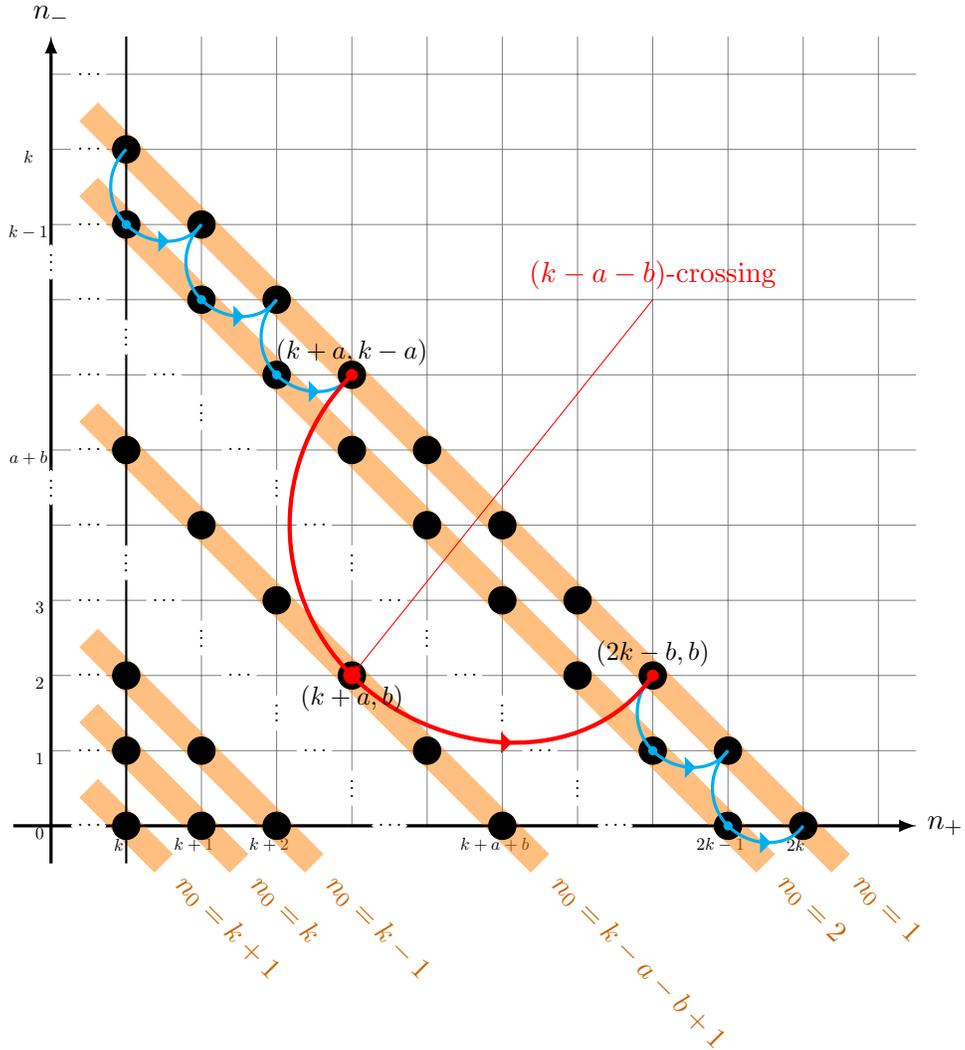
\begin{figure}
			\begin{center}
				\begin{tikzpicture}[scale=1,
									nullity_text/.style={right, 
														 rotate=-45, 
														 scale=1,
														 opacity=1, 
														 orange!80!black}, 
									nullity_line/.style={line width=3.5mm, 
														 orange, 
														 opacity=0.5},
									single_crossing_node/.style={draw, 
																  circle, 
																  fill,
																  scale = .3,
																  cyan},
									double_crossing_node/.style={draw, 
																  circle, 
																  fill,
																  scale = .6,
																  red}]
						\draw[help lines,step=1] (0,0) grid (11.5,10.5);
						\draw[-latex, very thick] (-.5,0) -- (11.5,0) node[right] {$n_+$};
						\draw[-latex, very thick] (0,-.5) -- (0,10.5) node[above] {$n_-$};
						\node[scale = .6] at (1-.1,-.25) {$k$};
						\node[scale = .6] at (2-.1,-.25) {$k+1$};
						\node[scale = .6] at (3-.1,-.25) {$k+2$};
						
						\node[scale = .6] at (6-.1,-.25) {$k+a+b$};
						
						\node[scale = .6] at (9-.1,-.25) {$2k-1$};
						\node[scale = .6] at (10-.1,-.25) {$2k$};
						
						\foreach \y in {0,1,2,3}
						{
							\node[scale = .6] at (-.15,\y-.1) {$\y$};
						}
						\node[scale = .6] at (-.3,5-.1) {$a+b$};
						
						\node[scale = .6] at (-.3,8-.1) {$k-1$};
						\node[scale = .6] at (-.3,9-.1) {$k$};
						
						\draw [nullity_line] (0.5,9.5) -- (10.5,-.5) node [nullity_text] {$n_0 = 1$};;		
						\draw [nullity_line] (0.5,8.5) -- (9.5,-.5) node [nullity_text] {$n_0 = 2$};;		
						
						\draw [nullity_line] (0.5,5.5) -- (6.5,-.5) node [nullity_text] {$n_0 = k-a-b+1$};;	

						\draw [nullity_line] (0.5,2.5) -- (3.5,-.5) node [nullity_text] {$n_0 = k-1$};;
						\draw [nullity_line] (0.5,1.5) -- (2.5,-.5) node [nullity_text] {$n_0 = k$};;	
						\draw [nullity_line] (0.5,.5) -- (1.5,-.5) node [nullity_text] {$n_0 = k+1$};;		
						\draw[-, thick] (-.5,0) -- (10.5,0);						
						\draw[-, thick] (1,-.5) -- (1,10.5);
						\foreach \h in {0,1,2, 5, 8,9}
						{
							\pgfmathsetmacro{\z}{9-\h}
							\foreach \y in {0,...,\h}
							{
								\node[draw, circle, fill] (start) at (1+\y,\h-\y) {};
							}
						}
						\pgfmathsetmacro{\h}{3.5}
						\pgfmathsetmacro{\z}{9-\h}
						\foreach \y in {-1,0,1,...,\h}
						{
							\node[draw, diamond,fill,white,scale=.9] () at (1+\y,\h-\y) {};
							\node[scale = .7] () at (1+\y,.07+\h-\y) {$\vdots$};
						}
						\foreach \y in {-.5,.5,...,\h}
						{
							\node[draw, diamond,fill,white,scale=.9] () at (1+\y,\h-\y) {};
							\node[scale = .7] () at (1.03+\y,\h-\y) {$\cdots$};
						}
						\foreach \y in {0,1,...,10}
						{
							\node[draw, diamond,fill,white,scale=.9] () at (.5,\y) {};
							\node[scale = .7] () at (.53,\y) {$\cdots$};
						}
						\pgfmathsetmacro{\h}{6.5}
						\pgfmathsetmacro{\z}{9-\h}
						\foreach \y in {-1,0,1,...,\h}
						{
							\node[draw, diamond,fill,white,scale=.9] () at (1+\y,\h-\y) {};
							\node[scale = .7] () at (1+\y,.07+\h-\y) {$\vdots$};
						}
						\foreach \y in {-.5,.5,...,\h}
						{
							\node[draw, diamond,fill,white,scale=.9] () at (1+\y,\h-\y) {};
							\node[scale = .7] () at (1.03+\y,\h-\y) {$\cdots$};
						}
					\foreach \x in {9,8,3,2,1}
					{
						\node[single_crossing_node] () at (\x,9-\x) {};
						\draw[-, very thick, cyan] (\x+1,9-\x) 
													to[out=235,in = -45] 
													node {\midarrow} 
													(\x,9-\x) 
													to[out=135,in = -135] 
													(\x,10-\x);
					}
					\node[double_crossing_node] () at (4,2) {};
					\draw[-, ultra thick, red] (8,2) 
												to[out=235,in = -45] 
												node {\midarrow} 
												(4,2) 
												to[out=135,in = -135] 
												(4,6);
					
					\draw[latex-, red] (4,2) -- (8,7) node[above] {$(k-a-b)$-crossing};

						\node[circle, draw, fill,red, scale = .4] () at (4,6) {};
						\node[scale=.9, black] () at (4,6.3) {$(k+a,k-a)$};
						
						\node[circle, draw, fill,red, scale = .4] () at (4,2) {};
						\node[scale=.9, black] () at (4,1.7) {$(k+a,b)$};
						
						\node[circle, draw, fill,red, scale = .4] () at (8,2) {};
						\node[scale=.9, black] () at (8,2.3) {$(2k-b,b)$};
						
				\end{tikzpicture}
			\end{center}
			\caption{The lattice of pairs $(n_+, n_-)$ of Laplacian inertias of $G \cdot G \cdot \cdots \cdot G$ where $G = K_3$ with two negative edges and a positive edge, and $G$ is repeated $k$ times.}
			\label{fig6}
			\end{figure}
			First note that the line $n_0 = 1$ has the equation $n_+ + n_- = 2k$. 	
			To illustrate the crossings giving all possible Laplacian inertias, start at the point $(k,k)$ on this line and go to the end point $(2k,0)$ passing through any other given point $(k+a,b)$ on the lattice, where $a$ and $b$ are nonnegative integers with $a+b < k$. From the point $(k,k)$, go through $a$ single crossings, then a $(k-a-b)$-crossing, and finally $b$ more single crossings to end at the point $(2k,0)$. In order to construct such sequence of crossings, we need to find weights $\gamma$ such that the crossing polynomial $M(\Gamma(t))$ of $\Gamma=(G,\gamma)$ has zeros
			\begin{equation*}
				t_1 < t_2 < \cdots < t_a < t_\ast < t_{a+1} < t_{a+2} < \cdots < t_{a+b},
			\end{equation*}
			where the $t_i$ have multiplicity one, and $t_\ast$ has multiplicity $k-a-b$. Note that, by Example \ref{ex_kay_three_crossing},
			\begin{equation*}
				M(r^-\Gamma(t)) = r^2 t (t - \frac{2}{r})
			\end{equation*}
			has a positive zero at $t = \frac{2}{r}$. Choose $r_i = \frac{2}{i}$, for $i = 1,2, \ldots, a+b$, and $r_\ast = \frac{4}{2a+1}$, and define $\Gamma_i = r_i^-\Gamma$, and $\Gamma_\ast = r_\ast^-\Gamma$. Then the crossing polynomial of 
			\begin{equation*}
				\Lambda = \Gamma_1 \cdot \Gamma_2 \cdot\, \cdots\, \cdot \Gamma_a \cdot \Gamma_\ast \cdot\, \cdots\, \cdot  \Gamma_\ast \cdot \Gamma_{a+1} \cdot\, \cdots\, \cdot \Gamma_{a+b},
			\end{equation*}
			where $\Gamma_\ast$ is repeated $k-a-b$ times, is 
			\begin{equation*}
				M(\Lambda(t)) = \frac{4^{2k-a-b}}{(2a+1)^{2(k-a-b)}} t^k 
				(t - \frac{2a+1}{2})^{k-a-b} \prod_{i=1}^{a+b} \frac{t - i}{i^2},
			\end{equation*}
			which has simple zeros at each $i = 1,2,\ldots,a+b$, and a zero of multiplicity $k-a-b$ at $\frac{2a+1}{2}$. 
Thus the Laplacian inertia corresponding to every point on the lattice can be achieved by choosing all non-negative pairs $a,b$ with $a + b < k$, and there are ${\tau +2 \choose 2}$ such points.
		\end{proof}

	\subsection{Graphs on a fixed number of vertices}
		\begin{lemma}
			The maximum flexibility of a signed graph on $n \geq 4$ vertices is $\tau = n - 1$ (i.e., if all eigenvalues cross).
		\end{lemma}
		\begin{proof}
			In order to maximize $\tau$ among all graphs on $n$ vertices, minimize $c_- + c_+$. This means the graph has a spanning tree with all negative edges and a spanning tree with all positive edges, which is possible for $n \geq 4$. In this case $\tau = n-1$.
		\end{proof}
		
		Note that the lower left-most point on the lattice of pairs $(n_+,n_-)$ of a graph corresponds to a $(\tau+1)$-crossing, i.e. when all $\tau+1$ negative eigenvalues cross at the same time. For a graph on $n$ vertices with $\tau = n-1$ this means that $n-1$ negative eigenvalues become zero simultaneously, and then all become positive. That is, at the instance of the crossing, $L(\Gamma(t))$ has all eigenvalues equal to $0$. But the only symmetric matrix with all eigenvalues equal to $0$ is the zero matrix. Thus the only graph $\Gamma$ on $n$ vertices that allows $\il(\Gamma) = (0,0,n)$ is the empty graph. This means the number of Laplacian inertias for connected signed graphs on $n \geq 2$ vertices is at most ${n+1 \choose 2} - 1$. 
		
		\begin{definition}
			Let $\sqcup$ denote the disjoint union binary operation on sets. Let $G_1 = (V_1, E_1)$ and $G_2 = (V_2, E_2)$ be signed graphs. Then the \emph{negative join} of two graphs $G_1 \vee^- G_2$ is a graph on the vertex set $V_1 \sqcup V_2$, the edge set $E_1 \sqcup E_2 \sqcup E_3$, where $E_3 = \{ uv \, | \, u \in V_1, v \in V_2  \}$, signs of the edges in $E_1$ and $E_2$ are the same as in $G_1$ and $G_2$, and the signs of the edges in $E_3$ are all negative. Informally, it is obtained by connecting with negative edges all the vertices in $G_1$ to all the vertices in $G_2$. 
		\end{definition}

		\begin{lemma}
			Let $G$ be a signed graph on $n$ vertices. Assume that there is a consistent weighting $\gamma$ of the edges of $G$ for which the Laplacian of $\Gamma = (G, \gamma)$ has rank $1$. Then up to isomorphism, $G$ is either $K_p \vee^- K_q$ or $- (K_p \vee^- K_q)$ for some $p,q \geq 1$ where $p+q = n$,  $\tau = n - 2$,
			and all edges of $K_p$ and $K_q$ are positive. 
\label{lemma3.10}		
		\end{lemma}	
		\begin{proof}
			Let $L$ be the Laplacian of $\Gamma = (G,\gamma)$ that has rank $1$.
			In a symmetric matrix with rank $1$ and no zero rows, all the rows are multiples of the 
			first row. Then all the diagonal entries of $L$ have the same sign, and no entries of $L$ 
			are zero, since a zero entry means there is an isolated vertex. First, assume that 
			$L_{1,1} > 0$ and that vertex $1$ is incident to $p-1$ positive edges and $q$ negative edges.
			Then the $p-1$ rows corresponding to the other vertices incident with these positive edges
			have the same sign pattern as row $1$, and the rest have the opposite sign pattern. 
			Hence, up to a permutation of rows and columns the matrix $L$ has the sign pattern
			\begin{equation*}
				\left[ 
					\begin{array}{c | c}
						\begin{array}{ccc}
							+ & \cdots & +	\\
							\vdots & \ddots & \vdots \\
							+ & \cdots & +						
						\end{array} & \begin{array}{ccc}
							- & \cdots & -	\\
							\vdots & \ddots & \vdots \\
							- & \cdots & -						
						\end{array} \\ \hline 
						\begin{array}{ccc}
							- & \cdots & -	\\
							\vdots & \ddots & \vdots \\
							- & \cdots & -						
						\end{array} & \begin{array}{ccc}
							+ & \cdots & +	\\
							\vdots & \ddots & \vdots \\
							+ & \cdots & +						
						\end{array}
					\end{array}
				\right],
			\end{equation*}
			where the top left block is $p \times p$ and the bottom right block is $q \times q$. 
	Hence $G = K_p \vee^- K_q $ with all edges of $K_p$ and $K_q$ positive, $c_+ = 2$ and $c_- = 1$.   		
			The case when $L_{1,1} < 0$ is similar and $c_+ = 1$ and $c_- = 2$. Finally, $\tau = n + 1 - c_- - c_+ = n - 2$.
		\end{proof}			

\noindent		
Note that the signed graphs $G$ in Lemma \ref{lemma3.10} are structurally balanced;
see, for example, \cite[Definition 2]{Shi}.
		
		\begin{example}	
		In this example we show that 
		$K_3 \vee^- K_3$ achieves all possible Laplacian inertias within the bounds of Theorem \ref{cor_inertia_bounds_disconnected}.	
		Consider the graph $G = K_3 \vee^- K_3$ where all edges of both graphs $K_3$ are positive (as shown in Figure \ref{fig7}(A)). Then $n = 6$, $c_+ = 2$, $c_- = 1$, and $\tau = 4$. Thus, the maximum number of possible Laplacian inertias for $G$ is ${\tau +2 \choose 2} = 15$, and all such inertias are shown in Figure \ref{fig7}(B). 
		\begin{figure}
		\begin{center}
			\begin{tikzpicture}[scale = 1,
								nullity_text/.style={right, 
													 rotate=-45, 
													 scale=1,
													 opacity=1, 
													 orange!80!black}, 
								nullity_line/.style={line width=3.5mm, 
													 orange, 
													 opacity=0.5},]
			\begin{scope}[shift={(-6,0)}]
			\node[] () at (-.5,5) {(A)};
				\node[draw,circle,scale=1] (1) at (2,0) {$1$};
				\node[draw,circle,scale=1] (2) at (2,2) {$2$};
				\node[draw,circle,scale=1] (3) at (2,4) {$3$};
				\draw[ultra thick, black] (1) to[out=45,in=-35] (2);
				\draw[ultra thick, black] (1) to[out=45,in=-35] (3);
				\draw[ultra thick, black] (2) to[out=45,in=-35] (3);
				\node[draw,circle,scale=1] (4) at (0,0) {$4$};
				\node[draw,circle,scale=1] (5) at (0,2) {$5$};
				\node[draw,circle,scale=1] (6) at (0,4) {$6$};
				\draw[ultra thick, black] (4) to[out=135,in=-135] (5);		
				\draw[ultra thick, black] (4) to[out=135,in=-135] (6);	
				\draw[ultra thick, black] (5) to[out=135,in=-135] (6);					
				\foreach \x in {1,2,3}
				{
					\foreach \y in {4,5,6}
					{
						\draw[dashed, ultra thick, red] (\x) -- (\y);
					}
				}
			\end{scope}										 
			\node[] () at (-1,5) {(B)};									 
				\draw[help lines,step=1] (0,0) grid (5.5,4.5);
				\draw[-latex, very thick] (-.5,0) -- (5.5,0) node[right] {$n_+$};
				\draw[-latex, very thick] (0,-.5) -- (0,4.5) node[above] {$n_-$};
			
				\foreach \x in {1,...,5}
				{
					\node[scale = .6] () at (\x-.1,-.25) {$\x$};
				}
				\foreach \x in {1,...,4}
				{
					\node[scale = .6] () at (-.25,\x-.1) {$\x$};
				}
				\draw[nullity_line] (.5,4.5) -- (5.5,-.5) node [nullity_text] {$n_0 = 1$};
				\draw[nullity_line] (.5,3.5) -- (4.5,-.5) node [nullity_text] {$n_0 = 2$};	
				\draw[nullity_line] (.5,2.5) -- (3.5,-.5) node [nullity_text] {$n_0 = 3$};
				\draw[nullity_line] (.5,1.5) -- (2.5,-.5) node [nullity_text] {$n_0 = 4$};
				\draw[nullity_line] (.5,.5) -- (1.5,-.5) node [nullity_text] {$n_0 = 5$};
				\foreach \x in {1,...,5}
				{
					\pgfmathsetmacro{\z}{5-\x}
					\foreach \y in {0,...,\z}
					{
						\node[draw, circle, fill] (start) at (\x,\y) {};
					}
				}

			\end{tikzpicture}
		\end{center}
		\caption{(A) The graph $K_3 \vee^- K_3$. (B) The lattice representing all possible Laplacian inertias for $K_3 \vee^- K_3$.}
		\label{fig7}
		\end{figure}
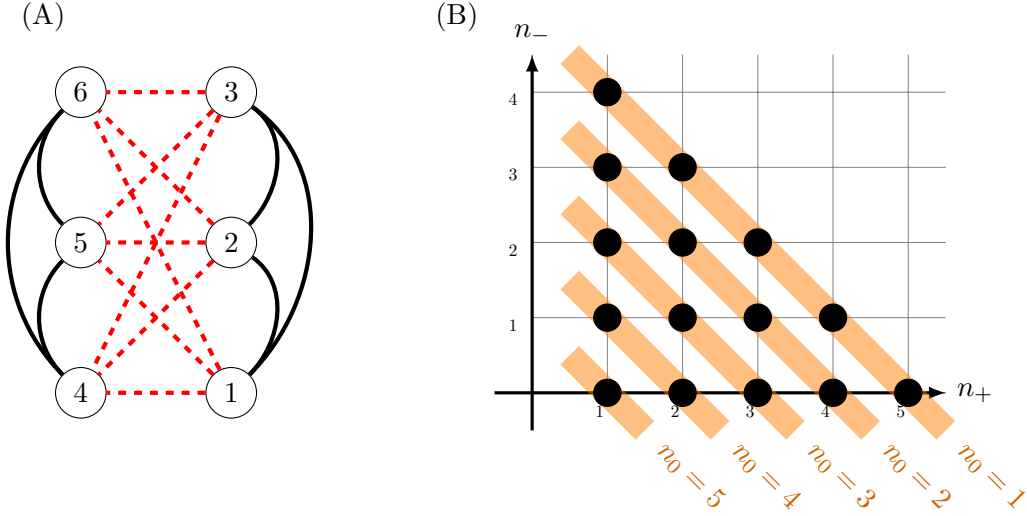
		First choose the weights of two of the positive edges in one $K_3$ and one of the positive edges on the other $K_3$ to be $2$, and choose all the other weights to be $\pm 1$ and consistent with the signs of the edges of $G$. Then 
		\begin{equation*}
			M(\Gamma(t)) = -9t (3t-4) (3t-5) (t-1) (t-2). 
		\end{equation*}
		That is, for $t = i/6$ with $i = 5,7, \ldots, 13$, $\Gamma(t)$ realizes all the Laplacian inertias with $n_0 = 1$, and for $t = i/3$ for $i = 3,4,5,6$, $\Gamma(t)$ realizes all the Laplacian inertias with $n_0 = 2$.
		
		Choose the weights of all of the edges of one $K_3$ to be $2$ and the weights on all the other edges to be $\pm 1$. Then 
		\begin{equation*}
			M(\Gamma(t)) = -81 t (t-1)^2 (t-2)^2.
		\end{equation*}
		This means for $t = 1$ and $t = 2$, $\Gamma(t)$ realizes Laplacian inertias $(1,2,3)$ and $(3,0,3)$, respectively. Now choose the weights on the edges of one $K_3$ to be $0.5$, $1$, and $2$, and the weights on all the other edges to be $\pm 1$. Then 
		\begin{equation*}
			M(\Gamma(t)) = \frac{-27}{2} t (6t^2 - 14t + 7) (t - 1)^2,
		\end{equation*}
		which has a zero of multiplicity $2$ at $t = 1$, and simple zeros at $t = \frac{7 \pm \sqrt{7}}{6} \approx 0.7257, 1.608$. That is, $\Gamma(1)$ realizes the Laplacian inertia $(2,1,3)$. This covers all the possible Laplacian inertias with nullity $3$.
		
		Now if $\gamma$ is chosen with all negative edges having weight $-1$, all but one positive edge having weight $1$, and the remaining positive edge having weight $2$ or $0.5$, then 
		\begin{equation*}
		 	M(\Gamma(t)) = -27t(3t-5)(t-1)^3
		\end{equation*}
		or
		\begin{equation*}
		 	M(\Gamma(t)) = -27t(3t-2)(t-1)^3,
		\end{equation*}
		which means for $t = 5/3$ (respectively $t = 2/3$), $\Gamma(t)$ realizes the Laplacian inertia $(1,1,4)$ (respectively $(2,0,4)$), which are all the possible Laplacian inertias with nullity $4$.
		
		Finally, choose a consistent weighting $\gamma$ of $G$ where all the weights are $\pm 1$. Then for $\Gamma = (G,\gamma)$, $\Gamma(1)$ realizes the last possible Laplacian inertia $(1,0,5)$ since
		\begin{equation*}
			M(\Gamma(t)) = -81t(t-1)^4.
		\end{equation*}
	\label{ex3.11}	
		\end{example}

		\begin{example}
		\label{ex_K4_max_tau}	
In this example we determine that 7 is the maximum number of Laplacian inertias achieved by a signed graph on 4 vertices.	
%
The maximum flexibility of a signed graph on 4 vertices is 3, and this is achieved if and only if $c_+ = c_- = 1$,
which implies that the signed graph is connected.
Up to isomorphism, 
the only signed graph on $4$ vertices with $c_\pm = 1$  is  shown in Figure \ref{fig8}(A).
For $n=4$,  the number of possible Laplacian inertias is at most ${5 \choose 2} - 1 = 9$. 
Since $\tau = 3$,  of these nine, $2 \tau + 1 = 7$  are realizable via single crossings. 
Since the graph is not isomorphic to $\pm (K_m \vee^-  K_n)$, it cannot achieve rank $1$, that is, it cannot obtain a $2$-crossing. Hence, only the seven Laplacian inertias mentioned above are realizable and their pairs $(n_+, n_-)$ are indicated by black dots in the lattice shown in Figure \ref{fig8}(B). 	
Any other signed graph on 4 vertices has $\tau \leq 2$, and hence at most $\binom{2+2}{2} = 6$ Laplacian 
inertias can be achieved.
Thus, the maximum number of Laplacian inertias that any signed graph on $4$ vertices can achieve is indeed $7$.

			\begin{figure}
			\begin{center}
				\begin{tikzpicture}[scale = 1,
									nullity_text/.style={right, 
														 rotate=-45, 
														 scale=1,
														 opacity=1, 
														 orange!80!black}, 
									nullity_line/.style={line width=3.5mm, 
														 orange, 
														 opacity=0.5},]
				\node[] () at (-1,4) {(B)};
						\draw[help lines,step=1] (0,0) grid (3.5,3.5);
						\draw[-latex, very thick] (-.5,0) -- (3.5,0) node[right] {$n_+$};
						\draw[-latex, very thick] (0,-.5) -- (0,3.5) node[above] {$n_-$};
					
						\foreach \x in {1,2,3}
						{
							\node[scale = .6] () at (\x-.1,-.25) {$\x$};
							\node[scale = .6] () at (-.25,\x-.1) {$\x$};
						}
						\draw[nullity_line] (-.5,3.5) -- (3.5,-.5) node [nullity_text] {$n_0 = 1$};
						\draw[nullity_line] (-.5,2.5) -- (2.5,-.5) node [nullity_text] {$n_0 = 2$};	
						\draw[nullity_line] (-.5,1.5) -- (1.5,-.5) node [nullity_text] {$n_0 = 3$};
						\draw[nullity_line] (-.5,0.5) -- (0.5,-.5) node [nullity_text] {$n_0 = 4$};
						\foreach \x in {0,...,3}
						{
							\pgfmathsetmacro{\z}{3-\x}
							\foreach \y in {0,...,\z}
							{
								\node[draw, circle, fill] (start) at (\x,\y) {};
							}
						}
						\foreach \x in {0,...,1}
						{
							\pgfmathsetmacro{\z}{1-\x}
							\foreach \y in {0,...,\z}
							{
								\node[draw, circle, fill, white] (start) at (\x,\y) {};
								\node[draw, circle, black] (start) at (\x,\y) {};
							}
						}

					\begin{scope}[shift={(-5,0)}]
					\node[] () at (-.5,4) {(A)};
						\node[draw, circle] (1) at (0,0) {$1$};
						\node[draw, circle] (2) at (0,2) {$2$};
						\node[draw, circle] (3) at (2,2) {$3$};
						\node[draw, circle] (4) at (2,0) {$4$};
						
						\draw[ultra thick, black] (1) -- (2) -- (3) -- (4);
						\draw[dashed, ultra thick, red] (2) -- (4) -- (1) -- (3);
					\end{scope}
				\end{tikzpicture}
			\end{center}
			\caption{(A) The only (up to isomorphism) signed graph $G$ on $4$ vertices with $c_\pm = 1$. (B) Black dots represent all achievable Laplacian inertias by $G$ and the white dots are the ones that cannot be achieved by a consistent weighting of $G$.}
			\label{fig8}	
			\end{figure}
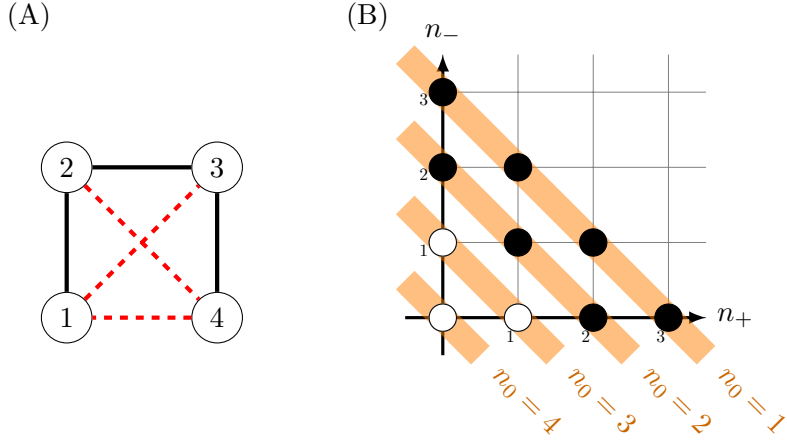
		
A signed graph on 4 vertices with $\tau = 2$ can achieve 6 Laplacian inertias if and only if it
allows rank 1: there must be a 2-crossing.
	The only connected signed graphs on $4$ vertices that could achieve rank $1$
			are the ones isomorphic to either $\pm(K_1 \vee^- K_3)$ or $\pm(K_2 \vee^- K_2)$. 
			 Five of the 6 possible Laplacian inertias are guaranteed by single crossings and, for example,  a $2$-crossing is given by the rank $1$ Laplacian matrix 			\begin{equation*}
				L = \left[
					\begin{array}{rrrr}
						9 &  -3 &  -3 &  -3 \\
						 -3 & 1 & 1 & 1 \\
						 -3 & 1 & 1 & 1 \\
						- 3 & 1 & 1 & 1 				
					\end{array}
				\right],
			\end{equation*}
which is a weighted signed Laplacian for the graph $K_1 \vee^- K_3$ with all edges of $K_3$ positive.

		For the signed graph $K_2 \vee^- K_2$, which has $n=4, c_+ = 2, c_- = 1$, and $\tau = 2$, reasoning as in Example \ref{ex3.11} shows that all ${\tau + 2 \choose 2} = 6$ possible Laplacian inertias are achieved.

		\end{example}
		
		In the following theorem we summarize the results of this section.
		
		\begin{theorem}
		Let $G$ be a signed graph on $n \geq 3$ vertices. 
		Then the number of  Laplacian inertias is at most ${n+1 \choose 2} - 3$.			
		\end{theorem}
		\begin{proof}
		Let $\tau$ be the flexibility of $G$.
			First note that for $n = 3$, the flexibility $\tau \leq 1$. Hence by Theorem \ref{thm_tau_bound} the statement is true. For $n \geq 4$ we prove the theorem in two cases. 
			\begin{itemize}
				\item If $\tau \leq n-2$, then the number of Laplacian inertias is at most ${\tau+2 \choose 2} \leq {n+1 \choose 2} - 3$.
				\item If $\tau = n-1$, then the number of Laplacian inertias is at most ${n+1 \choose 2} - 3$.
			\end{itemize}
			For a given signed graph $G$ on $n \geq 4$ vertices with $c_\pm$ connected components of $G_\pm$, Theorem \ref{cor_inertia_bounds_disconnected} provides sharp bounds for $n_\pm$. 
			Bronski et al. \cite{BD14,BDK15} show that realizations of a matrix with single crossings of the eigenvalues from the lower half plane to the upper half plane can always be found. There are $\tau = n + c - c_- - c_+$ single crossings. Starting from one extreme bound given by Corollary \ref{cor_inertia_bounds_disconnected}, for example with $t > 0$ sufficiently small, $t$ can be increased gradually to achieve $2 \tau + 1$ Laplacian inertias, $\tau$ of them with $n_0 = 2$, and $\tau + 1$ of them with $n_0 = 1$. 
		
			If $G$ is not isomorphic to $K_m \vee^- K_{n-m}$ for some integer $0 < m < n$, then $L(G, \gamma)$ cannot have nullities $0$ or $1$ for any consistent weighting $\gamma$. Hence it cannot achieve the ($n-1$)-crossing nor either of the two $(n-2)$-crossings. And if $G$ is isomorphic to $K_m \vee^- K_{n-m}$, then either $c_- = 1$ and $c_+ = 2$, or $c_- = 2$ and $c_+ = 1$. In both cases $\tau = n-2$, and since $n \geq 3$, by Theorem \ref{thm_tau_bound} it follows that ${\tau+2 \choose 2} = {n \choose 2} \leq {n+1 \choose 2} - 3$.
		\end{proof}
		
		We have shown that the maximum possible number of Laplacian inertias for a signed graph on $n$ vertices, namely ${n+1 \choose 2} - 3$, is achieved for $n = 3, 4$. But it is unknown whether this bound can be achieved for larger values of $n$. If a signed graph on $n$ vertices has a cut vertex, then it cannot achieve the maximum number of Laplacian inertias on $n$ vertices.
		
		For example, consider the graph $G$ in Example \ref{ex_K4_max_tau} depicted in Figure \ref{fig8}(A). The signed graph $G \cdot G$ has $7$ vertices, so ${n+1 \choose 2} - 3 = 25$. However, by Figure \ref{fig8}(B) and Remark \ref{cor_mink_sum_inertias}, the signed graph $G \cdot G$ achieves only the $18$ Laplacian inertias with $n_0 = 1, 2$ or $3$. Similarly, the signed graph $G \cdot G \cdot G$ has $10$ vertices and ${n+1 \choose 2} - 3 = 52$, but it achieves only the $34$ Laplacian inertias with nullities $n_0 = 1,2,3,4$. Note that realizing a Laplacian inertia with $n_0 \geq 5$ requires a $4$-crossing or a higher crossing, which means by the pigeon hole principle that $G$ must realize a $2$-crossing, but that is not possible. That is, $G \cdot G \cdot G$ cannot achieve any other Laplacian inertias.

\bigskip	
\noindent		
{\bf Acknowledgement}. We thank the referee for reading the manuscript and pointing out additional relevant references.

\bibliographystyle{unsrt}
\bibliography{refabv}

\bigskip
\begin{center}
k1monfared@gmail.com,\ \ \ gmacgill@math.uvic.ca,\ \ \  dolesky@uvic.ca,\ \ \  pvdd@math.uvic.ca
\end{center}
\end{document}